# Bifurcating Continued Fractions II

**Ashok Kumar Mittal**
**Department of Physics**
**Allahabad University, Allahabad – 211 002, India**
**(Email address: mittal_a@vsnl.com)**

**Ashok Kumar Gupta**
**Department of Electronics and Communication**
**Allahabad University, Allahabad - 211 002, INDIA**
**(Email address: akgjkiapt@hotmail.com)**

**Abstract:**

In an earlier paper[1] we introduced the notion of 'bifurcating continued fractions' in a heuristic manner. In this paper a formal theory is developed for the 'bifurcating continued fractions'.

Keywords: Continued fractions, Quadratic irrational, Cubic irrational, Golden Mean, Fibonacci sequence, Tribonacci sequence.

## 1. Introduction

Continued fractions provide deep insight into mathematical problems; particularly into the nature of numbers.[2] Continued fractions have found applications in various areas[3] of Physics such as Fabry-Perot interferometry, quasi-amorphous states of matter and chaos.

Any eventually periodic continued fraction represents a quadratic irrational. Conversely, Lagrange's theorem asserts that the continued fraction expansion of every quadratic irrational is eventually periodic. A purely periodic continued fraction represents a quadratic irrational of a special kind called a reduced quadratic irrational. A quadratic irrational is said to be reduced if it is greater than 1 and the other root of the quadratic equation that it satisfies, lies between -1 and 0. Conversely, the continued fraction expansion of a reduced quadratic irrational is purely periodic.

The continued fraction expansion consisting of the number 1 repeated indefinitely represents the 'golden mean'. This satisfies the quadratic equation $x^2 = x + 1$. The convergents of the continued fraction are obtained as the ratio of the successive terms of the Fibonacci sequence. Each term of the Fibonacci sequence is obtained by summing the previous two terms of the sequence.

A straightforward generalization of the Fibonacci sequence is that of the Tribonacci sequence, in which each term is obtained by summing the previous three terms of the sequence. Tribonacci number, the limiting ratio of the successive terms of the Tribonacci sequence satisfy the Tribonacci equation $x^3 = x^2 + x + 1$. However, the conventional continued fraction expansion of the Tribonacci number offers no satisfying pattern that would reflect the simple generalization from the 'golden mean'.

In our earlier paper[1], a method to generalize the continued fraction expansion so as to remove this deficiency of the conventional continued fractions was suggested in a heuristic manner. This method associates two coupled non-negative integer sequences with two real numbers and based on these sequences a 'bifurcating continued fraction' expansion is obtained. The 'bifurcating continued fraction' resembles the conventional continued fraction, except that both the numerator and the denominator bifurcate in a “Fibonacci Tree” like manner, whereas in conventional continued fractions only the denominator bifurcates. In particular, two sequences consisting of the number 1 repeated indefinitely represent the Tribonacci number. Other cubic variants[4] of the 'golden mean', like the Moore number satisfying the cubic equation $x^3 = x^2 + 1$, also find very simple and elegant representation as a 'bifurcating continued fraction'. Bifurcating continued fraction representations reveal the secret beauty of many numbers.

In this paper a formal theory of ‘Bifurcating Continued Fractions’ (BCF) is developed. In Sec 2 the BCF algorithm for obtaining two integer sequences with an ordered pair of positive real numbers $(\alpha, \beta)$ is stated. How the pair $(\alpha, \beta)$ can be recovered by use of a ‘Fibonacci Tree’ associated with the two integer sequences is described.

Sec 3 gives a formal development of the BCF theory. There are two processes, which are converse to each other. First, given an ordered pair of real numbers $(\alpha, \beta)$ one can obtain a pair of sequences of non-negative integers by the BCF algorithm. This pair of sequences is called the BCF expansion of $(\alpha, \beta)$. Second, given a pair of sequences of non-negative integers, one can associate a pair of Fibonacci Trees with it using simple bifurcating rules. These trees, when they are finite, have straightforward interpretation as a pair of rational numbers. If appropriate limits exist, the pair of Fibonacci Trees associated with an ordered pair of infinite sequences of non-negative integers may represent a pair of real numbers. In this case we say that the ordered pair of integer sequences is a Fibonacci Tree representation of $(\alpha, \beta)$. $(\alpha, \beta)$ may have more than one Fibonacci Tree representation.

In Theorem 1, it is shown that the BCF expansion of $(\alpha, \beta)$ is a Fibonacci Tree representation of $(\alpha, \beta)$. This provides formal justification for the heuristics of our earlier paper[1].

In order to establish Theorem 1, the notions of proper BCF representation and appropriate BCF representations were introduced. These are shown to be equivalent in Lemma 4. In Lemma 5, it is shown that the BCF expansion of $(\alpha, \beta)$ is the only proper BCF representation of $(\alpha, \beta)$.

In Lemma 6, necessary conditions for an ordered pair of non-negative integer sequences to be a proper BCF representation of an ordered pair of real numbers is obtained. In Theorem 2 these conditions are shown to be sufficient. By Lemma 5, if these conditions are satisfied by an ordered pair of non-negative integer sequences, they are the BCF expansions of $(\alpha, \beta)$ for some $\alpha$, $\beta$ and by Theorem 1 these sequences are the Fibonacci Tree representation of $(\alpha, \beta)$. As a corollary there is only one Fibonacci Tree representation of $(\alpha, \beta)$ which satisfies these conditions.

## 2. Generalizing the Continued Fractions

Given a positive real number $\alpha$, the continued fraction expansion of $\alpha$ is given by a sequence of non-negative integers $[a_0, a_1, \ldots\ldots. a_i \ldots.]$, obtained by the recurrence relation

$$a_i = \text{int}(\alpha_i),$$
$$\alpha_{i+1} = 1/(\alpha_i - a_i)\,, \qquad (1)$$

where $\alpha_0 = \alpha$.

Our generalization of the method of continued fraction is based on a generalization of equation (1). In this generalization we obtain an ordered pair of integer sequences $\{\{a_0, a_1,\ldots,a_i, \ldots\}, \{b_0, b_1,\ldots..,b_i,\ldots\}\}$ from a given ordered pair $\{\alpha,\beta\}$ of positive real numbers by using the recurrence relations

$$a_i = \text{int}(\alpha_i),$$
$$b_i = \text{int}(\beta_i),$$
$$\alpha_{i+1} = 1/(\beta_i - b_i),$$
$$\beta_{i+1} = (\alpha_i - a_i)/\beta_i - b_i), \qquad (2)$$

where $\alpha_0 = \alpha$ and $\beta_0 = \beta$.

The numbers $\alpha$ and $\beta$ can be recovered heuristically from the two integer sequences $\{a_i\}$ and $\{b_i\}$ by using the recurrence relations

$$\alpha_i = a_i + (\beta_{i+1} / \alpha_{i+1}),$$
$$\beta_i = b_i + (1/\alpha_{i+1}) \qquad (3)$$

Given the integer sequences $\{a_i\}$ and $\{b_i\}$, one can write down the bifurcating continued fractions by the bifurcating rule:

$$a_i \longrightarrow a_i < \begin{matrix} b_{i+1} \\ a_{i+1} \end{matrix}$$

and

$$b_i \longrightarrow b_i < \begin{matrix} 1 \\ a_{i+1} \end{matrix}$$

where

$$p < \begin{matrix} q \\ r \end{matrix} \quad = \quad p + \frac{q}{r}$$

To get $\alpha$, one begins with $a_0$ and to get $\beta$ one begins with $b_0$. In this notation

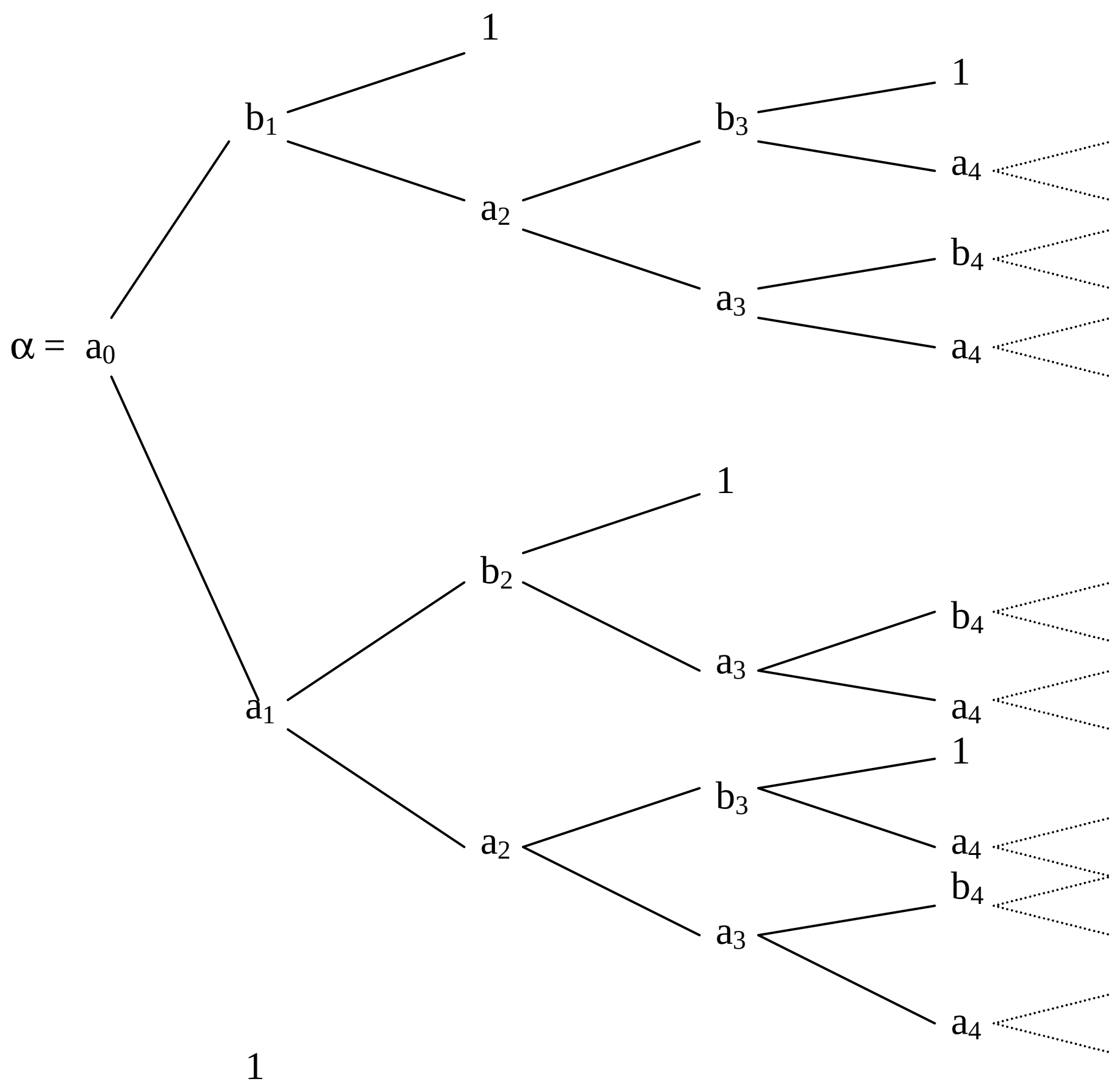

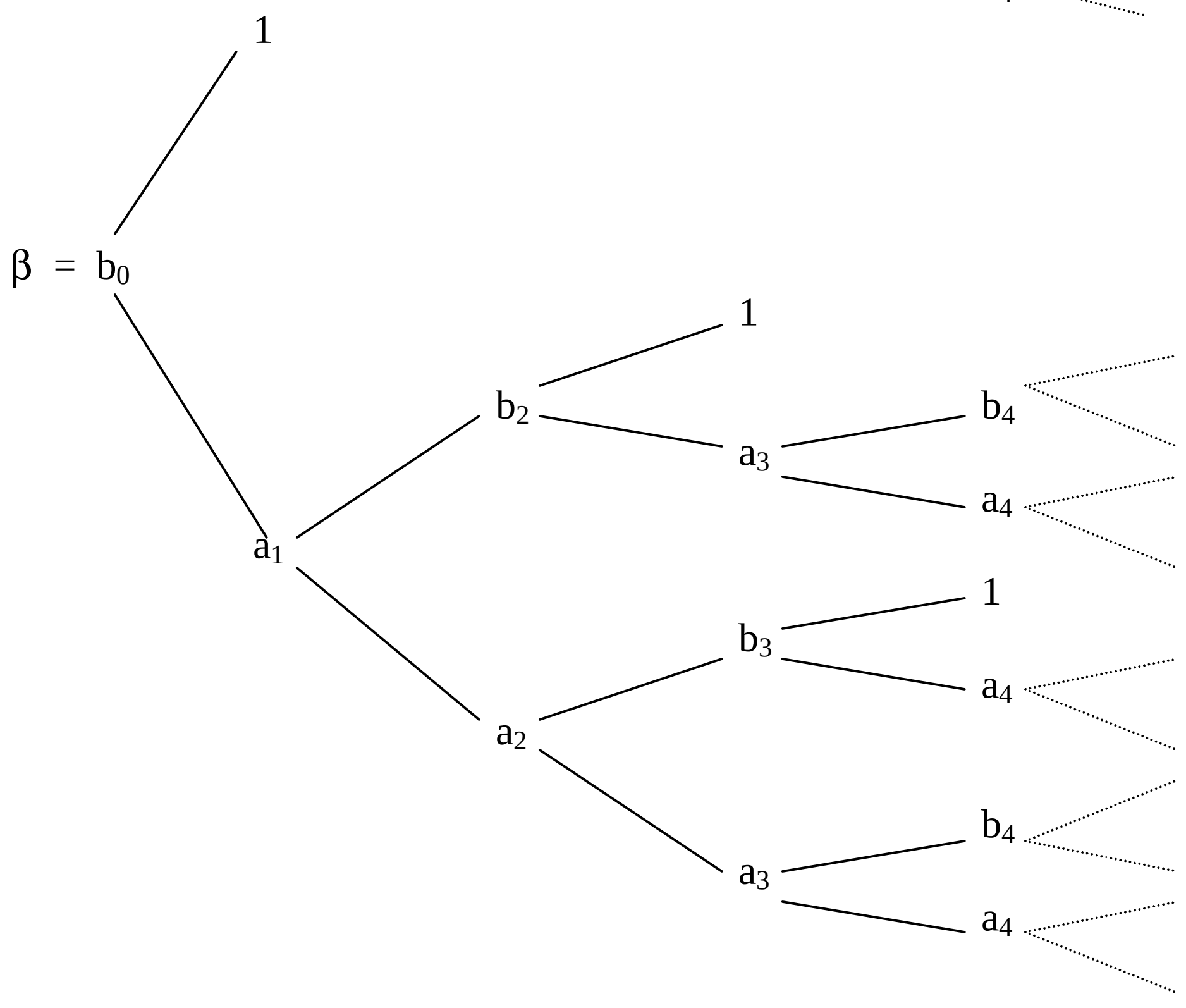

The Fibonacci tree structure of these expansions is obvious, in as much as at any branching level, the number of a's (#a) equals the sum of number of a's at the previous two levels. The same holds for the b's.

$$\#a_n = \#a_{n-1} + \#a_{n-2} ,$$

$$\#b_n = \#b_{n-1} + \#b_{n-2}.$$

## 3. Formal development of Bifurcating Continued Fraction Theory

In our earlier paper[1] we demonstrated how the BCF algorithm represented by the set of equations (2) appear to be a very attractive candidate for generalizing the conventional continued fraction algorithm. Hence there is a need to develop a formal theory of BCF along lines parallel to that of the conventional continued fractions.

**Def 1**. Given an ordered pair $\{\alpha, \beta\}$ of positive real numbers, the ordered pair of integer sequences $\{\{a_0, a_1,\ldots,a_i, \ldots\}, \{b_0, b_1,\ldots,b_i,\ldots\}\}$ obtained by the Bifurcating Continued Fraction algorithm of equations (2) is called the BCF expansion of $\{\alpha, \beta\}$. The algorithm is undefined if $\beta_i$ becomes an integer for some i = n (say); in this case we say that the BCF expansion of $\{\alpha, \beta\}$ is finite and is given by the ordered pair of finite sequences $\{\{a_0, a_1,\ldots,a_i, \ldots a_{n-1}, \alpha_n\}, \{b_0, b_1,\ldots,b_i,\ldots b_n\}\}$ where $a_i$, $b_i$ are integers and $\alpha_n$ is real.

**Lemma 1**. If $\{\{a_0, a_1,\ldots,a_i, \ldots\}, \{b_0, b_1,\ldots,b_i,\ldots\}\}$ is the BCF expansion of the ordered pair $\{\alpha, \beta\}$, then $1 \leq a_i \geq b_i$ for $i \geq 1$ and if $a_n = b_n$ for some $n \geq 1$ then $b_{n+1} = 0$.

Proof: We have,

$$\begin{aligned} a_i &= \text{int}(\alpha_i), \\ b_i &= \text{int}(\beta_i), \\ \alpha_{i+1} &= 1/(\beta_i - b_i), \\ \beta_{i+1} &= (\alpha_i - a_i)/(\beta_i - b_i) , \end{aligned} \qquad (4)$$

where $\alpha_0 = \alpha$ and $\beta_0 = \beta$. Hence,

$$\alpha_{i+1} > 1 \quad \text{for } i \geq 0 \qquad (5)$$

Therefore,

$$a_i = \text{int}(\alpha_i) \geq 1 \quad \text{for } i \geq 1 \qquad (6)$$

Further,

$$\alpha_{i+1}/\beta_{i+1} = 1/(\alpha_i - a_i) > 1 \quad \text{for } i \geq 0 \tag{7}$$

Therefore,

$$\alpha_i > \beta_i \quad \text{for } i \geq 1 \tag{8}$$

Hence,

$$\text{int}(\alpha_i) \geq \text{int}(\beta_i) \quad \text{i.e.,} \quad a_i \geq b_i \quad \text{for } i \geq 1 \tag{9}$$

If $a_n = b_n$ for some $n \geq 1$, then by (4) and (8), we must have $\beta_{n+1} > 1$ and therefore $b_{n+1} \neq 0$.

**Lemma 2** If $\alpha$, $\beta$ are positive rational numbers, the BCF expansion of the ordered pair $\{\alpha, \beta\}$ is finite.

Proof: $\alpha$, $\beta$ are positive rational numbers implies $\alpha_i$, $\beta_i$ are positive rational for all i. Hence one can find integers $u_i$, $v_i$, $w_i$ such that

$$\alpha_i = u_i/w_i \quad \text{and} \quad \beta_i = v_i/w_i \tag{10}$$

Moreover, $1 < \alpha_i > \beta_i$ for $i \geq 1$. Therefore,

$$w_i < u_i > v_i \quad \text{for } i \geq 1 \tag{11}$$

Then,

$$\begin{aligned} u_{i+1}/w_{i+1} &= \alpha_{i+1} = 1/(\beta_i - b_i) = 1/\{(v_i/w_i) - b_i\} = w_i/(v_i - b_i w_i) \\ v_{i+1}/w_{i+1} &= \beta_{i+1} = (\alpha_i - a_i)/(\beta_i - b_i) = \{(u_i/w_i) - a_i\}/\{(v_i/w_i) - b_i\} = (u_i - a_i w_i)/(v_i - b_i w_i) \end{aligned} \tag{12}$$

These are satisfied by

$$u_{i+1} = w_i, \quad v_{i+1} = u_i - a_i w_i, \quad w_{i+1} = v_i - b_i w_i \tag{13}$$

It follows from the BCF algorithm that $w_{i+1} < w_i$ for all i. Thus $w_i$ is a strictly decreasing sequence of natural numbers and therefore there must exist an n such that $w_{n+1} = 0$. But that implies $\beta_n = v_n/w_n = b_n$ is an integer and therefore the BCF algorithm terminates.

Def 2. The Fibonacci Tree sum of an ordered pair of sequences $\{\{a_0, a_1,\ldots,a_n, \alpha_{n+1}\}, \{b_0, b_1,\ldots,b_n, \beta_{n+1}\}\}$ where $a_i$, $b_i$ are non-negative integers and $\alpha_{n+1}$, $\beta_{n+1}$ are positive reals, is defined as the ordered pair of real numbers $(\alpha, \beta) = (\alpha_0, \beta_0)$ denoted by $[\{a_0, a_1,\ldots,a_n, \alpha_{n+1}\}, \{b_0, b_1,\ldots, b_n, \beta_{n+1}\}]$ ) and obtained by the recurrence relations

$[\{a_0, a_1,\ldots,a_n, \alpha_{n+1}\}, \{b_0, b_1,\ldots,b_n, \beta_{n+1}\}]$

$$= [\{a_0, a_1,\ldots,a_{n-1}, a_n + (\beta_{n+1}/\alpha_{n+1})\}, \{b_0, b_1,\ldots,b_{n-1}, b_n + (1/\alpha_{n+1}\}] \quad (14)$$

with the terminating condition

$$[\{\alpha_0\}, \{\beta_0\}] = (\alpha_0, \beta_0) \quad (15)$$

Lemma 3. If $[\{a_0, a_1,\ldots,a_n, \alpha_{n+1}\}, \{b_0, b_1,\ldots,b_n, \beta_{n+1}\}] = [\{a_0, a_1,\ldots,a_n, \alpha'_{n+1}\}, \{b_0, b_1,\ldots,b_n, \beta'_{n+1}\}]$ then $\alpha'_{n+1} = \alpha_{n+1}$ and $\beta'_{n+1} = \beta_{n+1}$.

Proof: The Lemma holds for n = 0. For, if $[\{a_0, \alpha_1\}, \{b_0, \beta_1\}] = [\{a_0, \alpha'_1\}, \{b_0, \beta'_1\}]$ then

$$a_0 + (\beta_1/\alpha_1) = a_0 + (\beta'_1/\alpha'_1) \quad (16)$$

and $$b_0 + (1/\alpha_1) = b_0 + (1/\alpha'_1) \quad (17)$$

whence equation (17) implies $\alpha'_1 = \alpha_1$ and then equation (16) yields $\beta'_1 = \beta_1$.

Assume as an induction hypothesis that the Lemma holds for n = 0, 1, ….. , k-1, so that if $[\{a_0, a_1,\ldots,a_{k-1}, \alpha_k\}, \{b_0, b_1,\ldots,b_{k-1}, \beta_k\}] = [\{a_0, a_1,\ldots,a_{k-1}, \alpha'_k\}, \{b_0, b_1,\ldots,b_{k-1}, \beta'_k\}]$ then $\alpha'_k = \alpha_k$ and $\beta'_k = \beta_k$.

If $[\{a_0, a_1,\ldots,a_k, \alpha_{k+1}\}, \{b_0, b_1,\ldots,b_k, \beta_{k+1}\}] = [\{a_0, a_1,\ldots,a_k, \alpha'_{k+1}\}, \{b_0, b_1,\ldots,b_k, \beta'_{k+1}\}]$ then

$[\{a_0, a_1,\ldots,a_{k-1}, a_k + (\beta_{k+1}/\alpha_{k+1})\}, \{b_0, b_1,\ldots,b_{k-1}, b_k + (1/\alpha_{k+1}\}]$

$$= [\{a_0, a_1,\ldots,a_{k-1}, a_k + (\beta'_{k+1}/\alpha'_{k+1})\}, \{b_0, b_1,\ldots,b_{k-1}, b_k + (1/\alpha'_{k+1}\}] \quad (18)$$

By induction hypothesis eqn (18) implies

$$a_k + (\beta_{k+1}/\alpha_{k+1}) = a_k + (\beta'_{k+1}/\alpha'_{k+1}) \quad (19)$$

and $$b_k + (1/\alpha_{k+1}) = b_k + (1/\alpha'_{k+1}) \quad (20)$$

whence equation (20) implies $\alpha'_{k+1} = \alpha_{k+1}$ and then equation (19) yields $\beta'_{k+1} = \beta_{k+1}$.

Thus the Lemma holds for n = 0 and further, if the Lemma holds for n = 0, 1, ….., k-1, it holds for n = k. Hence the Lemma holds for all n.

Def 3. The n-term Fibonacci Tree sum of an ordered pair of non-negative integer sequences ($\{a_0, a_1,\ldots,a_i, \ldots \}, \{b_0, b_1,\ldots,b_i,\ldots \}$) is defined by

$$(\alpha^{(n)}, \beta^{(n)}) = [\{a_0, a_1,\ldots,a_{n-1}, a_n\}, \{b_0, b_1,\ldots,b_{n-1}, b_n\}] \quad (21)$$

If $\lim_{n\to\infty} (\alpha^{(n)}, \beta^{(n)}) = (\alpha, \beta)$ the ordered pair of integer sequences is called a Fibonacci Tree representation of $(\alpha, \beta)$.

Example: The ordered pair of sequences ($\{1, 1, 1, \ldots \}, \{2, 2, 2, \ldots \}$) is a Fibonacci Tree representation of $(\alpha, \beta)$ where $\alpha$ satisfies $\alpha^3 = \alpha^2 + 2\alpha + 1$ and $\beta = 2 + (1/\alpha)$. However the BCF expansion of $(\alpha, \beta)$ is ($\{2, 2, 3, 2, 3, 2, 3, 2, \ldots \}, \{2, 0, 0, 0, \ldots \}$).

Def 4. An ordered pair of non-negative integer sequences ($\{a_0, a_1,\ldots,a_i, \ldots \}, \{b_0, b_1, \ldots , b_i,\ldots \}$) is called a proper BCF representation of $(\alpha, \beta)$, where $\alpha, \beta$ are positive real numbers, if for every n there exist positive reals $\alpha_n$ and $\beta_n$ such that

$$(\alpha, \beta) = [\{a_0, a_1,\ldots,a_{n-1}, \alpha_n\}, \{b_0, b_1,\ldots,b_{n-1}, \beta_n\}] \quad (22)$$

and $1 < \alpha_n > \beta_n$ for all $n \geq 1$.

Def 5. An ordered pair of non-negative integer sequences ($\{a_0, a_1,\ldots,a_i, \ldots \}, \{b_0, b_1, \ldots,b_i,\ldots \}$) is called an appropriate BCF representation of $(\alpha, \beta)$, where $\alpha, \beta$ are positive real numbers, if for every n there exist positive reals $\alpha_n$ and $\beta_n$ such that

$$(\alpha, \beta) = [\{a_0, a_1,\ldots,a_{n-1}, \alpha_n\}, \{b_0, b_1,\ldots,b_{n-1}, \beta_n\}] \quad (23)$$

and $\text{int}(\alpha_n) = a_n$ , $\text{int}(\beta_n) = b_n$.

Lemma 4. An ordered pair of non-negative integer sequences ($\{a_0, a_1,\ldots,a_i, \ldots \}, \{b_0, b_1,\ldots,b_i,\ldots \}$) is a proper BCF representation of $(\alpha, \beta)$ iff it is an appropriate BCF representation of $(\alpha, \beta)$

Proof: Suppose ($\{a_0, a_1,\ldots,a_i, \ldots \}, \{b_0, b_1,\ldots,b_i,\ldots \}$) is a proper BCF representation of $(\alpha, \beta)$. Then for every n there exist positive reals $\alpha_n$ and $\beta_n$ such that

$$(\alpha, \beta) = [\{a_0, a_1,\ldots,a_{n-1}, \alpha_n\}, \{b_0, b_1,\ldots,b_{n-1}, \beta_n\}] \quad (24)$$

and $1 < \alpha_n > \beta_n$ for all $n \geq 1$.

Clearly, we must have

$$[\{a_0, a_1,\ldots,a_n, \alpha_{n+1}\}, \{b_0, b_1,\ldots,b_n, \beta_{n+1}\}] = [\{a_0, a_1,\ldots,a_{n-1}, \alpha_n\}, \{b_0, b_1,\ldots,b_{n-1}, \beta_n\}] \quad (25)$$

Applying eqn (14) and Lemma 3 to eqn (25) we get,

$$\alpha_n = a_n + (\beta_{n+1} / \alpha_{n+1}), \quad \beta_n = b_n + (1/\alpha_{n+1}) \tag{26}$$

Since $0 < \beta_{n+1}/\alpha_{n+1} < 1$ and $0 < 1/\alpha_{n+1} < 1$ for all n, it follows from eqns (26) that $a_n = \text{int}(\alpha_n)$ and $b_n = \text{int}(\beta_n)$. So the given pair of non-negative integer sequences is an appropriate BCF representation of $(\alpha, \beta)$.

Conversely, suppose that the given pair of non-negative integer sequences is an appropriate BCF representation of $(\alpha, \beta)$. Then there exist reals $\alpha_n$ and $\beta_n$ satisfying eqns (24) and (26). Moreover, $a_n = \text{int}(\alpha_n)$ and $b_n = \text{int}(\beta_n)$. But then from eqns (26) it follows that $1 < \alpha_n > \beta_n$ for all $n \geq 1$. So the given pair of non-negative integer sequences is a proper BCF representation of $(\alpha, \beta)$.

Lemma 5. The BCF expansion of $(\alpha, \beta)$ is the only proper BCF representation of $(\alpha, \beta)$.

Proof: The BCF expansion of $(\alpha, \beta)$ gives rise to real numbers $\alpha_i$, $\beta_i$ by the algorithm

$$\alpha_0 = \alpha, \ \beta_0 = \beta$$

$$a_i = \text{int}(\alpha_i), \ b_i = \text{int}(\beta_i) \tag{27}$$

$$\alpha_{i+1} = 1/(\beta_i - b_i), \ \beta_{i+1} = (\alpha_i - a_i)/(\beta_i - b_i)$$

It is straightforward to verify that for every n

$$(\alpha, \beta) = [\{a_0, a_1, \ldots, a_{n-1}, \alpha_n\}, \{b_0, b_1, \ldots, b_{n-1}, \beta_n\}] \tag{28}$$

and $$\alpha_{n+1} = 1/(\beta_n - b_n) > 1, \ \beta_{n+1}/\alpha_{n+1} = (\alpha_n - a_n) < 1 \tag{29}$$

so that $1 < \alpha_n > \beta_n$ for all $n \geq 1$. Hence the BCF expansion of $(\alpha, \beta)$ is a proper BCF representation of $(\alpha, \beta)$.

Let $\{\{a'_0, a'_1, \ldots, a'_i, \ldots\}, \{b'_0, b'_1, \ldots, b'_i, \ldots\}\}$ be another proper BCF representation of $(\alpha, \beta)$. Then for every n there exist positive reals $\alpha'_n$ and $\beta'_n$ such that

$$(\alpha, \beta) = [\{a'_0, a'_1, \ldots, a'_{n-1}, \alpha'_n\}, \{b'_0, b'_1, \ldots, b'_{n-1}, \beta'_n\}] \tag{30}$$

and $1 < \alpha'_n > \beta'_n$ for all $n \geq 1$. Then by eqn (26)

$$\alpha'_n = a'_n + (\beta'_{n+1} / \alpha'_{n+1}), \quad \beta'_n = b'_n + (1/\alpha'_{n+1}) \tag{31}$$

But

$$\alpha'_0 = \alpha_0 = a_0 + (\beta_1/\alpha_1), \;\; \beta'_0 = \beta_0 = b_0 + (1/\alpha_1) \tag{32}$$

Hence

$$a_0 + (\beta_1/\alpha_1) = a'_0 + (\beta'_1/\alpha'_1), \;\; b_0 + (1/\alpha_1) = b'_0 + (1/\alpha'_1) \tag{33}$$

Because $b_0$, $b'_0$ are integers and $\alpha_1, \alpha'_1 > 1$ it follows that $b_0 = b'_0$ and $\alpha_1 = \alpha'_1$. Again because $a_0$, $a'_0$ are integers and $\beta_1/\alpha_1$ and $\beta'_1/\alpha'_1 < 1$ it follows that $a_0 = a'_0$ and $\beta_1 = \beta'_1$. By induction it is straightforward to show that $a'_n = a_n$ and $b_n = b'_n$ for all n.

Lemma 6. If an ordered pair of non-negative integer sequences $(\{a_0, a_1,\ldots,a_i, \ldots\}, \{b_0,b_1,\ldots,b_i,\ldots\})$ is a proper BCF representation of $(\alpha, \beta)$, then (i) $1 \leq a_n \geq b_n$ (ii) $a_n = b_n$ for any $n \geq 1$ implies $b_{n+1} \neq 0$.

Proof: For every n there exist positive reals $\alpha_n$ and $\beta_n$ such that

$$(\alpha, \beta) = [\{a_0, a_1,\ldots,a_{n-1}, \alpha_n\}, \{b_0, b_1,\ldots,b_{n-1}, \beta_n\}] \tag{34}$$

and $1 < \alpha_n > \beta_n$ if $n \geq 1$.

(i) By Lemma 4,

$$a_n = \text{int}(\alpha_n) \geq 1$$

and $$a_n = \text{int}(\alpha_n) \geq \text{int}(\beta_n) = b_n$$

(ii) Suppose $a_n = b_n$ for some $n \geq 1$ and $b_{n+1} = 0$. Then,

$$\beta_{n+1} = (b_{n+1} + 1/\alpha_{n+2}) < 1$$

so that $\alpha_n = (a_n + \beta_{n+1}/\alpha_{n+1}) < (b_n + 1/\alpha_{n+1}) = \beta_n$

which contradicts the hypothesis.

Lemma 7. $[\{a_m, a_{m+1},\ldots,a_n\}, \{b_m, b_{m+1},\ldots,b_n\}] = ((A_{m,n}/A_{m+1,n}),(B_{m,n}/A_{m+1,n}))$ where $A_{m,n}$ and $B_{m,n}$ satisfy the recurrence relations $A_{m,n} = a_m A_{m+1,n} + b_{m+1} A_{m+2,n} + A_{m+3,n}$ and $B_{m,n} = b_m A_{m+1,n} + A_{m+2,n}$, $0 \leq m \leq n < \infty$, with the notional initial conditions, $A_{n+3,n} = 0$, $A_{n+2,n} = 0$, $A_{n+1,n} = 1$.

Proof: The Lemma holds for m = n, because $[\{a_n\}, \{b_n\}] = (a_n, b_n) = ((A_{n,n}/A_{n+1,n}),(B_{n,n}/A_{n+1,n}))$. The Lemma also holds for m = n-1, because

$$[\{a_{n-1}, a_n\}, \{b_{n-1}, b_n\}] = (a_{n-1} + (b_n/a_n), b_{n-1} + (1/a_n)) = ((a_{n-1}a_n + b_n)/a_n, (b_{n-1}a_n + 1)/a_n)$$

$$= ((a_{n-1} A_{n,n} + b_n A_{n+1,n} + A_{n+2,n})/A_{n,n}, (b_{n-1} A_{n,n} + A_{n+1,n})/A_{n,n})$$

$$= ((A_{n-1,n}/A_{n,n}),(B_{n-1,n}/ A_{n,n})) \qquad (35)$$

Suppose that the Lemma holds for $m = n-1, n-2,\ldots\ldots,n-k$. Then,

$$[\{a_{n-k}, a_{n-k+1},\ldots,a_n\}, \{b_{n-k}, b_{n-k+1},\ldots..,b_n\}] = ((A_{n-k,n}/A_{n-k+1,n}),(B_{n-k,n}/ A_{n-k+1,n})) \qquad (36)$$

where

$$A_{n-k,n} = a_{n-k} A_{n-k+1,n} + b_{n-k+1} A_{n-k+2,n} + A_{n-k+3,n} \qquad (37)$$

and $\quad B_{n-k,n} = b_{n-k} A_{n-k+1,n} + A_{n-k+2,n} \qquad (38)$

Therefore,

$$[\{a_{n-k-1}, a_{n-k},\ldots,a_n\}, \{b_{n-k-1}, b_{n-k},\ldots..,b_n\}]$$

$$= [\{a_{n-k-1}, A_{n-k,n}/A_{n-k+1,n}\}, \{b_{n-k-1}, B_{n-k,n}/ A_{n-k+1,n}\}]$$

$$= ( a_{n-k-1} + (B_{n-k,n}/ A_{n-k+1,n})/ (A_{n-k,n}/A_{n-k+1,n}), b_{n-k-1} + (1/ (A_{n-k,n}/A_{n-k+1,n}))$$

$$= ( a_{n-k-1} + (B_{n-k,n}/ A_{n-k,n}), b_{n-k-1} + (A_{n-k+1,n} / A_{n-k,n}))$$

$$= (( a_{n-k-1} A_{n-k,n} + b_{n-k} A_{n-k+1,n} + A_{n-k+2,n})/ A_{n-k,n}, (b_{n-k-1} A_{n-k,n} + A_{n-k+1,n})/ A_{n-k,n})$$

$$= ((A_{n-k-1,n}/A_{n-k,n}),(B_{n-k-1,n}/ A_{n-k,n})) \qquad (39)$$

where

$$A_{n-k-1,n} = a_{n-k-1} A_{n-k,n} + b_{n-k} A_{n-k+1,n} + A_{n-k+2,n} \qquad (40)$$

and $\quad B_{n-k-1,n} = b_{n-k-1} A_{n-k,n} + A_{n-k+1,n} \qquad (41)$

Hence the Lemma holds for $m = n - k-1$ and therefore, by induction for all non-negative integer values of m less than or equal to n.

Lemma 8. $A_{m,n} = a_n A_{m,n-1} + b_n A_{m,n-2} + A_{m,n-3}$

Proof: The Lemma holds for $n = m$ because by Lemma 7 and the notional initial conditions we have,

$$A_{m,m} = a_m A_{m+1,m} + b_{m+1} A_{m+2,m} + A_{m+3,m} = a_m$$

$$= a_m A_{m,m-1} + b_m A_{m,m-2} + A_{m,m-3} \qquad (42)$$

The Lemma holds for n = m+1 because

$$A_{m,m+1} = a_m A_{m+1,m+1} + b_{m+1} A_{m+2,m+1} + A_{m+3,m+1}$$

$$= a_m a_{m+1} + b_{m+1} = a_{m+1} a_m + b_{m+1}$$

$$= a_{m+1} A_{m,m} + b_{m+1} A_{m,m-1} + A_{m,m-2} \qquad (43)$$

The Lemma holds for n = m + 2 because

$$A_{m,m+2} = a_m A_{m+1,m+2} + b_{m+1} A_{m+2,m+2} + A_{m+3,m+2}$$

$$= a_m (a_{m+1} a_{m+2} + b_{m+2}) + b_{m+1} a_{m+2} + 1$$

$$= a_{m+2} (a_m a_{m+1} + b_{m+1}) + b_{m+2} a_m + 1$$

$$= a_{m+2} A_{m,m+1} + b_{m+2} A_{m,m} + A_{m,m-1} \qquad (44)$$

Assume as an induction hypothesis that the Lemma holds for n = m, m+1, m+2, …, m+k so that,

$$A_{m,m+k} = a_{m+k} A_{m,m+k-1} + b_{m+k} A_{m,m+k-2} + A_{m,m+k-3} \qquad (45)$$

We have by Lemma 7,

$$A_{m,m+k+1} = a_m A_{m+1,m+k+1} + b_{m+1} A_{m+2,m+k+1} + A_{m+3,m+k+1}$$

$$\begin{aligned} &= a_m (a_{m+k+1} A_{m+1,m+k} + b_{m+k+1} A_{m+1,m+k-1} + A_{m+1,m+k-2}) \\ &\quad + b_{m+1} (a_{m+k+1} A_{m+2,m+k} + b_{m+k+1} A_{m+2,m+k-1} + A_{m+2,m+k-2}) \\ &\quad + (a_{m+k+1} A_{m+3,m+k} + b_{m+k+1} A_{m+3,m+k-1} + A_{m+3,m+k-2}) \end{aligned}$$

$$\begin{aligned} &= a_{m+k+1} (a_m A_{m+1,m+k} + b_{m+1} A_{m+2,m+k} + A_{m+3,m+k}) \\ &\quad + b_{m+k+1} (a_m A_{m+1,m+k-1} + b_{m+1} A_{m+2,m+k-1} + A_{m+3,m+k-1}) \\ &\quad + (a_m A_{m+1,m+k-2} + b_{m+1} A_{m+2,m+k-2} + A_{m+3,m+k-2}) \end{aligned}$$

$$= a_{m+k+1} A_{m,m+k} + b_{m+k+1} A_{m,m+k-1} + A_{m,m+k-2} \qquad (46)$$

Hence the Lemma holds for n = m + k + 1 and therefore for all integer values of n greater than or equal to m.

Lemma 9. $B_{m,n} = a_n B_{m,n-1} + b_n B_{m,n-2} + B_{m,n-3}$

Proof: We have by Lemma 7 and 8

$$\begin{aligned} B_{m,n} &= b_m A_{m+1,n} + A_{m+2,n} \\ &= b_m (a_n A_{m+1,n-1} + b_n A_{m+1,n-2} + A_{m+1,n-3}) + (a_n A_{m+2,n-1} + b_n A_{m+2,n-2} + A_{m+2,n-3}) \\ &= a_n (b_m A_{m+1,n-1} + A_{m+2,n-1}) + b_n (b_m A_{m+1,n-2} + A_{m+2,n-2}) + (b_m A_{m+1,n-3} + A_{m+2,n-3}) \\ &= a_n B_{m,n-1} + b_m B_{m,n-2} + B_{m,n-3} \qquad (47) \end{aligned}$$

Lemma 10.
$$\det \begin{vmatrix} A_n & A_{n-1} & A_{n-2} \\ B_n & B_{n-1} & B_{n-2} \\ C_n & C_{n-1} & C_{n-2} \end{vmatrix} = 1$$

where $A_i = A_{0,i}$ , $B_i = B_{0,i}$ and $C_i = A_{1,i}$

Proof: For n > 2, using Lemmas 8 and 9, we have

$$\det \begin{vmatrix} A_n & A_{n-1} & A_{n-2} \\ B_n & B_{n-1} & B_{n-2} \\ C_n & C_{n-1} & C_{n-2} \end{vmatrix} = \det \begin{vmatrix} a_n A_{n-1} + b_n A_{n-2} + A_{n-3} & A_{n-1} & A_{n-2} \\ a_n B_{n-1} + b_n B_{n-2} + B_{n-3} & B_{n-1} & B_{n-2} \\ a_n C_{n-1} + b_n C_{n-2} + C_{n-3} & C_{n-1} & C_{n-2} \end{vmatrix}$$

$$= \det \begin{vmatrix} A_{n-3} & A_{n-1} & A_{n-2} \\ B_{n-3} & B_{n-1} & B_{n-2} \\ C_{n-3} & C_{n-1} & C_{n-2} \end{vmatrix} = \det \begin{vmatrix} A_{n-1} & A_{n-2} & A_{n-3} \\ B_{n-1} & B_{n-2} & B_{n-3} \\ C_{n-1} & C_{n-2} & C_{n-3} \end{vmatrix} = \ldots\ldots = \ldots\ldots$$

$$= \det \begin{vmatrix} A_2 & A_1 & A_0 \\ B_2 & B_1 & B_0 \\ C_2 & C_1 & C_0 \end{vmatrix} = \det \begin{vmatrix} A_{0,2} & A_{0,1} & A_{0,0} \\ B_{0,2} & B_{0,1} & B_{0,0} \\ A_{1,2} & A_{1,1} & A_{1,0} \end{vmatrix}$$

$$= \det \begin{vmatrix} a_0 A_{1,2} + b_1 A_{2,2} + 1 & a_0 A_{1,1} + b_1 & a_0 \\ b_0 A_{1,2} + A_{2,2} & b_0 A_{1,1} + 1 & b_0 \\ A_{1,2} & A_{1,1} & 1 \end{vmatrix} \qquad \text{(using Lemma 7)}$$

$$= \det \begin{vmatrix} b_1 A_{2,2} + 1 & b_1 & 0 \\ A_{2,2} & 1 & 0 \\ A_{1,2} & A_{1,1} & 1 \end{vmatrix} = \det \begin{vmatrix} 1 & 0 & 0 \\ A_{2,2} & 1 & 0 \\ A_{1,2} & A_{1,1} & 1 \end{vmatrix} = 1 \qquad (48)$$

Cor 1: $A_n$, $B_n$, $C_n$ cannot have a common factor

Cor 2: $(\alpha^{(n)} - \alpha^{(n-1)})(\beta^{(n-1)} - \beta^{(n-2)}) - (\alpha^{(n-1)} - \alpha^{(n-2)})(\beta^{(n)} - \beta^{(n-1)}) = 1/(C_{n-2}C_{n-1}C_n)$ (49)

Proof: Using def 3, Lemma 7 and Lemma 9, we have

$$\det \begin{vmatrix} A_n & A_{n-1} & A_{n-2} \\ B_n & B_{n-1} & B_{n-2} \\ C_n & C_{n-1} & C_{n-2} \end{vmatrix} = C_nC_{n-1}C_{n-2} \begin{vmatrix} \alpha^{(n)} & \alpha^{(n-1)} & \alpha^{(n-2)} \\ \beta^{(n)} & \beta^{(n-1)} & \beta^{(n-2)} \\ 1 & 1 & 1 \end{vmatrix} = 1 \qquad (50)$$

Therefore,

$$\alpha^{(n)}(\beta^{(n-1)} - \beta^{(n-2)}) + \alpha^{(n-1)}(\beta^{(n-2)} - \beta^{(n)}) + \alpha^{(n-2)}(\beta^{(n)} - \beta^{(n-1)}) = 1/(C_nC_{n-1}C_{n-2}) \qquad (51)$$

which yields eqn (49).

## Lemma 11.

If $1 \leq a_i \geq b_i$ for all $i \geq 1$ then the ordered pair of integer sequences $(\{a_0, a_1, \ldots, a_i, \ldots\}, \{b_0, b_1, \ldots, b_i, \ldots\})$ is a Fibonacci Tree representation of $(\alpha, \beta)$ for some real numbers $\alpha$ and $\beta$.

Proof: In view of definition 3, Lemmas 7 - 9 and notation of Lemma 1, we have to show that

$$\lim_{n\to\infty} \alpha^{(n)} \equiv \lim_{n\to\infty} A_n/C_n \qquad (52)$$

and

$$\lim_{n\to\infty} \beta^{(n)} \equiv \lim_{n\to\infty} B_n/C_n \qquad (53)$$

exist, where

$$\begin{aligned} A_n &= a_n A_{n-1} + b_n A_{n-2} + A_{n-3} \\ B_n &= a_n B_{n-1} + b_n B_{n-2} + B_{n-3} \\ C_n &= a_n C_{n-1} + b_n C_{n-2} + C_{n-3} \end{aligned} \qquad (54)$$

From eqns (54), we get,

$$\begin{aligned} (A_n/A_{n-1} - C_n/C_{n-1}) &= b_n (A_{n-2}/A_{n-1} - C_{n-2}/C_{n-1}) + (A_{n-3}/A_{n-1} - C_{n-3}/C_{n-1}) \\ &= b_n (C_{n-1}/C_{n-2} - A_{n-1}/A_{n-2}) (A_{n-2}/A_{n-1})( C_{n-2}/C_{n-1} ) \\ &\quad + (C_{n-1}/C_{n-3} - A_{n-1}/A_{n-3}) (A_{n-3}/A_{n-1} )(C_{n-3}/C_{n-1}) \end{aligned} \qquad (55)$$

But from eqns (54), we also get,

$$(A_n/A_{n-2} - C_n/C_{n-2}) = a_n (A_{n-1}/A_{n-2} - C_{n-1}/C_{n-2}) + (A_{n-3}/A_{n-2} - C_{n-3}/C_{n-2}) \qquad (56)$$

Using (56) in (55), we get,

$$(A_n/A_{n-1} - C_n/C_{n-1}) = b_n (C_{n-1}/C_{n-2} - A_{n-1}/A_{n-2})(A_{n-2}/A_{n-1})(C_{n-2}/C_{n-1}) + \{a_{n-1}(C_{n-2}/C_{n-3} - A_{n-2}/A_{n-3}) + (C_{n-4}/C_{n-3} - A_{n-4}/A_{n-3})\}(A_{n-3}/A_{n-1})(C_{n-3}/C_{n-1})$$

$$= b_n (C_{n-1}/C_{n-2} - A_{n-1}/A_{n-2})(A_{n-2}/A_{n-1})(C_{n-2}/C_{n-1}) + a_{n-1}(C_{n-2}/C_{n-3} - A_{n-2}/A_{n-3}) (A_{n-3}/A_{n-1})(C_{n-3}/C_{n-1}) + (A_{n-3}/A_{n-4} - C_{n-3}/C_{n-4})(A_{n-4}/A_{n-1})(C_{n-4}/C_{n-1}) \quad (57)$$

Therefore,

$$(A_n/C_n - A_{n-1}/C_{n-1})(C_n/A_{n-1}) = b_n (A_{n-2}/C_{n-2} - A_{n-1}/C_{n-1})(C_{n-2}/A_{n-1}) + a_{n-1}(A_{n-3}/C_{n-3} - A_{n-2}/C_{n-2})(C_{n-2}/A_{n-1})(C_{n-3}/C_{n-1}) + (A_{n-3}/C_{n-3} - A_{n-4}/C_{n-4})(C_{n-3}/A_{n-1})(C_{n-4}/C_{n-1}) \quad (58)$$

yielding

$$(A_n/C_n - A_{n-1}/C_{n-1}) = b_n (A_{n-2}/C_{n-2} - A_{n-1}/C_{n-1})(C_{n-2}/C_n) + a_{n-1}(A_{n-3}/C_{n-3} - A_{n-2}/C_{n-2})(C_{n-2}/C_n)(C_{n-3}/C_{n-1}) + (A_{n-3}/C_{n-3} - A_{n-4}/C_{n-4})(C_{n-3}/C_n)(C_{n-4}/C_{n-1}) \quad (59)$$

Now,

$$b_n C_{n-2}/C_n = b_n C_{n-2} / (a_nC_{n-1} + b_n C_{n-2} + C_{n-3}) = 0 \text{ ( if } b_n = 0)$$

$$= 1/\{(a_n/b_n)(C_{n-1}/C_{n-2}) + 1 + (1/b_n) (C_{n-3}/C_{n-2})\} \text{ (if } b_n \neq 0)$$

$$< 1/2 \qquad \text{(because, } C_{n-1} > C_{n-2}) \quad (60)$$

and

$$a_{n-1}(C_{n-2}/C_n)(C_{n-3}/C_{n-1}) = (C_{n-3}/C_n) \{a_{n-1} C_{n-2}/(a_{n-1}C_{n-2} + b_{n-1} C_{n-3} + C_{n-4})$$

$$< C_{n-3}/C_n \quad (61)$$

But $C_{n-3}/C_n = C_{n-3} / (a_nC_{n-1} + b_n C_{n-2} + C_{n-3})$

$$= C_{n-3} / \{a_n(a_{n-1}C_{n-2} + b_{n-1} C_{n-3} + C_{n-4}) + b_n (a_{n-2}C_{n-3} + b_{n-2} C_{n-4} + C_{n-5}) + C_{n-3}\}$$

$$< C_{n-3}/\{a_na_{n-1}C_{n-2} + (b_na_{n-2} + b_{n-1} + 1)C_{n-3}\}$$

$$< 1/(a_na_{n-1} + b_na_{n-2} + b_{n-1} + 1) \qquad \text{(because } a_n \geq 1, C_{n-2} > C_{n-3})$$

$$< \begin{cases} 1/4 & \text{if } b_n \neq 0, b_{n-1} \neq 0 \\ 1/3 & \text{if } b_n = 0, b_{n-1} \neq 0 \text{ or } b_n \neq 0, b_{n-1} = 0 \\ 1/2 & \text{if } b_n = b_{n-1} = 0 \end{cases} \qquad (62)$$

Similarly,

$$C_{n-4}/C_{n-1} < \begin{cases} 1/4 & \text{if } b_{n-1} \neq 0, b_{n-2} \neq 0 \\ 1/3 & \text{if } b_{n-1} = 0, b_{n-2} \neq 0 \text{ or } b_{n-1} \neq 0, b_{n-2} = 0 \\ 1/2 & \text{if } b_{n-1} = b_{n-2} = 0 \end{cases} \qquad (63)$$

Combining (62) and (63)

$$(C_{n-3}/C_n)(C_{n-4}/C_{n-1}) < \begin{cases} 1/16 & \text{if } b_n \neq 0, b_{n-1} \neq 0, b_{n-2} \neq 0 \\ 1/12 & \text{if } b_n \neq 0, b_{n-1} \neq 0, b_{n-2} = 0 \text{ or } b_n = 0, b_{n-1} \neq 0, b_{n-2} \neq 0 \\ 1/9 & \text{if } b_n = 0, b_{n-1} \neq 0, b_{n-2} = 0 \text{ or } b_n \neq 0, b_{n-1} = 0, b_{n-2} \neq 0 \\ 1/6 & \text{if } b_n \neq 0, b_{n-1} = 0, b_{n-2} = 0 \text{ or } b_n = 0, b_{n-1} = 0, b_{n-2} \neq 0 \\ 1/4 & \text{if } b_n = 0, b_{n-1} = 0, b_{n-2} = 0 \end{cases} \qquad (64)$$

Let,

$$\Delta_n = |(A_n/C_n - A_{n-1}/C_{n-1})| \qquad (65)$$

Then, by (60), (61), (62) and (64), we get,

$$\Delta_n < \begin{cases} (1/2)\,\Delta_{n-2} + (1/4)\,\Delta_{n-3} & (\text{if } b_n = 0) \\ (1/2)\,\Delta_{n-1} + (1/4)\,\Delta_{n-2} + (1/12)\,\Delta_{n-3} & (\text{if } b_n \neq 0, b_{n-1} \neq 0) \\ (1/2)\,\Delta_{n-1} + (1/3)\,\Delta_{n-2} + (1/9)\,\Delta_{n-3} & (\text{if } b_n \neq 0, b_{n-1} \neq 0, b_{n-2} \neq 0) \\ (1/2)\,\Delta_{n-1} + (1/3)\,\Delta_{n-2} + (1/6)\,\Delta_{n-3} & (\text{if } b_n \neq 0, b_{n-1} = 0, b_{n-2} = 0) \end{cases} \qquad (66)$$

Therefore,

$$\Delta_n < \begin{cases} D_n & \text{if } b_n \neq 0, b_{n-1} = 0, b_{n-2} = 0 \\ (17/18)\,D_n & \text{otherwise} \end{cases} \qquad (67)$$

where

$$D_n = \max\{\Delta_{n-1}, \Delta_{n-2}, \Delta_{n-3}\} \qquad (68)$$

Then,

$$D_{n+1} = \max\{\Delta_n, \Delta_{n-1}, \Delta_{n-2}\} \leq D_n \qquad (69)$$

Hence $\{D_n\}$ is a non-increasing sequence of positive numbers. We show that there exists $r < 1$ such that $D_{n+4} < r D_n$, wherefrom it follows that $\{D_n\}$ converges to zero. There are two cases to be considered. If $b_n \neq 0$, $b_{n-1} = 0$, $b_{n-2} = 0$ by (66), we get,

$$\Delta_{n+1} < (1/2)\,\Delta_n + (1/4)\,\Delta_{n-1} + (1/12)\,\Delta_{n-2}$$

$$< (5/6)\,D_n \qquad \text{(by (67) and (68))} \qquad (70)$$

$$\Delta_{n+2} < (1/2)\,\Delta_{n+1} + (1/3)\,\Delta_n + (1/9)\,\Delta_{n-1}$$

$$< (31/36)\,D_n \qquad \text{(by (67),(68) and (70))} \qquad (71)$$

$$\Delta_{n+3} < (1/2)\,\Delta_{n+2} + (1/3)\,\Delta_{n+1} + (1/6)\,\Delta_n$$

$$< (7/8)\,D_n \qquad \text{(by (67), (70) and (71))} \qquad (72)$$

Thus,

$$D_{n+4} = \max\{\Delta_{n+3}, \Delta_{n+2}, \Delta_{n+1}\} < (7/8)\,D_n \qquad (73)$$

Otherwise, by (67),

$$\Delta_n < (17/18)\,D_n \qquad (74)$$

$$\Delta_{n+1} < (1/2)\,\Delta_n + (1/3)\,\Delta_{n-1} + (1/6)\,\Delta_{n-2}$$

$$< (35/36)\,D_n \qquad \text{(by (68) and (74))} \qquad (75)$$

$$\Delta_{n+2} < (1/2)\,\Delta_{n+1} + (1/3)\,\Delta_n + (1/6)\,\Delta_{n-1}$$

$$< (209/216)\,D_n \qquad \text{(by (68), (74) and (75))} \qquad (76)$$

Therefore,

$$D_{n+3} = \max\{\Delta_{n+2}, \Delta_{n+1}, \Delta_n\} < (35/36)\,D_n \qquad (77)$$

Hence, by (69),

$$D_{n+4} \leq D_{n+3} < (35/36)\,D_n \qquad (78)$$

Thus $\{D_n\}$ is a non-increasing sequence with $D_{n+4} < r D_n$ where $r = 35/36 < 1$. Hence $\{D_n\}$ converges to zero as $n \to \infty$. From (71) it follows that the sequence $\{\Delta_n\}$ converges to zero, whence by eqn (65), $\lim_{n\to\infty} \alpha^{(n)} \equiv \lim_{n\to\infty} A_n/C_n$ exists. Replacing $A_n$ by $B_n$ in the above proof, it can be shown that $\lim_{n\to\infty} \beta^{(n)} \equiv \lim_{n\to\infty} B_n/C_n$ exists. Hence if $1 \leq a_i \geq b_i$ for all $i \geq 1$, then the ordered pair of integer sequences $(\{a_0, a_1, \ldots, a_i, \ldots\}, \{b_0, b_1, \ldots, b_i, \ldots\})$ is a Fibonacci Tree representation of $(\alpha, \beta)$ for some real numbers $\alpha$ and $\beta$.

Theorem 1. The BCF expansion of an ordered pair of real numbers $(\alpha, \beta)$ is a Fibonacci Tree representation of $(\alpha, \beta)$.

Proof: Let $(\{a_0, a_1, \ldots, a_i, \ldots\}, \{b_0, b_1, \ldots, b_i, \ldots\})$ be the BCF expansion of the ordered pair $(\alpha, \beta)$. Then by Lemma 1, $1 \leq a_i \geq b_i$ for all $i \geq 1$ and therefore by Lemma 11, there exist real numbers $\alpha'$ and $\beta'$ such that

$$\lim_{n\to\infty} \alpha^{(n)} = \alpha' \tag{79}$$

and
$$\lim_{n\to\infty} \beta^{(n)} = \beta' \tag{80}$$

where $\alpha^{(n)} = (A_n / C_n)$ and $\beta^{(n)} = (B_n / C_n)$

We have to show that $\alpha' = \alpha$ and $\beta' = \beta$.

Now,

$$\begin{aligned}|\alpha - \alpha^{(n)}| &= |(\alpha_n A_{n-1} + \beta_n A_{n-2} + A_{n-3})/(\alpha_n C_{n-1} + \beta_n C_{n-2} + C_{n-3}) \\ &\qquad - (a_n A_{n-1} + b_n A_{n-2} + A_{n-3})/(a_n C_{n-1} + b_n C_{n-2} + C_{n-3})| \\ &= |(A_n + A^{\sim}_n)/(C_n + C^{\sim}_n) - (A_n/C_n)|\end{aligned} \tag{81}$$

where

$$A^{\sim}_n = (\alpha_n - a_n) A_{n-1} + (\beta_n - b_n) A_{n-2} \tag{82}$$

and
$$C^{\sim}_n = (\alpha_n - a_n) C_{n-1} + (\beta_n - b_n) C_{n-2} \tag{83}$$

Then,

$$\begin{aligned}|\alpha - \alpha^{(n)}| &= |(A^{\sim}_n C_n - A_n C^{\sim}_n)/\{C_n (C_n + C^{\sim}_n)\}| \\ &< (1/C_n^2) |A^{\sim}_n C_n - A_n C^{\sim}_n| \\ &= (A_n / C_n) |(A^{\sim}_n / A_n) - (C^{\sim}_n / C_n)| \\ &= (A_n / C_n) |(\alpha_n - a_n)\{(A_{n-1}/A_n) - (C_{n-1}/C_n) + (\beta_n - b_n)\{(A_{n-2}/A_n) - (C_{n-2}/C_n) \\ &\leq \alpha^{(n)} \{|(\alpha^{(n-1)}/\alpha^{(n)} - 1|(C_{n-1}/C_n) + |(\alpha^{(n-2)}/\alpha^{(n)} - 1|(C_{n-2}/C_n)\}\end{aligned} \tag{84}$$

Hence,

$$|\alpha - \alpha'| = \lim_{n\to\infty} |\alpha - \alpha^{(n)}| = 0 \quad \text{(because } \lim \alpha^{(n)} = \alpha'\text{)} \tag{85}$$

Therefore $\alpha' = \alpha$. Similarly it can be shown that $\beta' = \beta$.

Theorem 2. If $1 \le a_i \ge b_i$ for all $i \ge 1$ and if $a_n = b_n$ implies $b_{n+1} \ne 0$ then the ordered pair of integer sequences $(\{a_0, a_1, \ldots, a_i, \ldots\}, \{b_0, b_1, \ldots, b_i, \ldots\})$ is a proper BCF representation of $(\alpha, \beta)$ for some real numbers $\alpha$ and $\beta$. By Lemma 5, this must be the BCF expansion of $(\alpha, \beta)$. By Theorem 1, this must be a Fibonacci Tree representation of $(\alpha, \beta)$

Proof: By Lemma11, there exists real numbers $\alpha$ and $\beta$ such that

$$(\alpha, \beta) = \lim_{n\to\infty} (\alpha^{(n)}, \beta^{(n)})$$

$$= \lim_{n\to\infty} [\{a_0, a_1, \ldots, a_n\}, \{b_0, b_1, \ldots, b_n\}] \tag{86}$$

From the proof of Lemma 11, it is clear that for every m, $\lim_{n\to\infty} [\{a_m, a_{m+1}, \ldots, a_n\}, \{b_m, b_{m+1}, \ldots, b_n\}]$ exists. Let this limit be denoted by $(\alpha_m, \beta_m)$. Thus for every m there exists $\alpha_m, \beta_m$ such that

$$(\alpha, \beta) = [\{a_0, a_1, \ldots, a_{m-1}, \alpha_m\}, \{b_0, b_1, \ldots, b_{m-1}, \beta_m\}] \tag{87}$$

We have to show that $1 < \alpha_m > \beta_m$. Since,

$$\alpha_m = \lim_{n\to\infty} (A_{m,n} / A_{m+1,n}) \tag{88}$$

and

$$\beta_m = \lim_{n\to\infty} (B_{m,n} / A_{m+1,n}) \tag{89}$$

we show that

$$(A_{m,n} / A_{m+1,n}) > 1 \tag{90}$$

and

$$(A_{m,n} / B_{m,n}) > 1 \tag{91}$$

so that $1 < \alpha_m > \beta_m$ (92)

By Lemma 7,

$$A_{m,n} = a_m A_{m+1,n} + b_{m+1} A_{m+2,n} + A_{m+3,n} \tag{93}$$

where given m, $A_{m+3,n} > 1$ for sufficiently large n.

Thus,

$$(A_{m,n} / A_{m+1,n}) = a_m + b_{m+1}(A_{m+2,n}/A_{m+1,n}) + (A_{m+3,n}/A_{m+1,n}) > 1 \qquad (94)$$

Again by Lemma 7,

$$(A_{m,n} / B_{m,n}) = (a_m A_{m+1,n} + b_{m+1} A_{m+2,n} + A_{m+3,n})/(b_m A_{m+1,n} + A_{m+2,n}) \qquad (95)$$

If $a_m > b_m$,

$$(A_{m,n} / B_{m,n}) > 1 \quad \text{(because } A_{m+1,n} > A_{m+2,n} \text{ by (90))} \qquad (96)$$

If $a_m = b_m$, then $b_{m+1} \neq 0$ and again

$$(A_{m,n} / B_{m,n}) > 1 \qquad (97)$$

Cor: There is only one Fibonacci Tree representation of $(\alpha, \beta)$ which satisfies the conditions (i) $1 \leq a_i \geq b_i$ for all $i \geq 1$ and (ii) $a_n = b_n$ implies $b_{n+1} \neq 0$

Theorem 3. If $(\alpha, \beta)$ has a purely periodic BCF expansion, then $\alpha$ satisfies a polynomial of degree at most three.

Proof: Let the BCF expansion have period n + 1. Then,

$$(\alpha, \beta) = [\{a_0, a_1, \ldots\ldots, a_n, \alpha\}, \{b_0, b_1, \ldots\ldots, b_n, \beta\}] \qquad (98)$$

Hence,

$$\alpha = (\alpha A_n + \beta A_{n-1} + A_{n-2})/(\alpha C_n + \beta C_{n-1} + C_{n-2}) \qquad (99)$$

$$\beta = (\alpha B_n + \beta B_{n-1} + B_{n-2})/ (\alpha C_n + \beta C_{n-1} + C_{n-2}) \qquad (100)$$

From (99) we get,

$$C_n\alpha^2 + C_{n-1}\alpha\beta + (C_{n-2} - A_n)\alpha - A_{n-1}\beta - A_{n-2} = 0 \qquad (101)$$

From (100) we get,

$$C_{n-1}\beta^2 + C_n\alpha\beta - B_n\alpha + (C_{n-2} - B_{n-1})\beta - B_{n-2} = 0 \qquad (102)$$

From (101) we get,

$$\beta = \{C_n\alpha^2 + (C_{n-2} - A_n)\alpha - A_{n-2}\}/(A_{n-1} - C_{n-1}\alpha) \qquad (103)$$

Putting (103) in (102),

$$C_{n-1}\{C_n\alpha^2 + (C_{n-2} - A_n)\alpha - A_{n-2}\}^2 + (C_n\alpha + C_{n-2} - B_{n-1})\{C_n\alpha^2 + (C_{n-2} - A_n)\alpha - A_{n-2}\}(A_{n-1} - C_{n-1}\alpha) - (B_n\alpha + B_{n-2})(A_{n-1} - C_{n-1}\alpha)^2 = 0 \tag{104}$$

The coefficient of $\alpha^4$ vanishes, therefore $\alpha$ satisfies an equation that is at most of degree three.

**Lemma 12.** $\alpha_{n+1}(\alpha_n \quad \beta_n \quad 1) = (\alpha_{n+1} \quad \beta_{n+1} \quad 1)\begin{pmatrix} a_n & b_n & 1 \\ 1 & 0 & 0 \\ 0 & 1 & 0 \end{pmatrix}$

Proof: R. H. S. $= ((\alpha_{n+1} a_n + \beta_{n+1}) \quad (\alpha_{n+1} b_n + 1) \quad \alpha_{n+1})$

$= \alpha_{n+1}((a_n + \beta_{n+1}/\alpha_{n+1}) \quad (b_n + 1/\alpha_{n+1}) \quad 1)$

$= \alpha_{n+1}(\alpha_n \quad \beta_n \quad 1)$

$=$ L. H. S.

Cor 1: $\alpha_{n+1}\mu_n^T = \mu_{n+1}^T R_n$

where $\mu_n = (\alpha_n \quad \beta_n \quad 1)^T$, $(\alpha_n, \beta_n) = [\{a_n, a_{n+1}, \ldots..\},\{b_n, b_{n+1}, \ldots..\}]$,

$$R_n = \begin{pmatrix} a_n & b_n & 1 \\ 1 & 0 & 0 \\ 0 & 1 & 0 \end{pmatrix} \tag{105}$$

Cor 2: $\mu_0^T = (\mu_1^T R_0)/\alpha_1 = (\mu_2^T R_1 R_0)/(\alpha_2\alpha_1) = \ldots\ldots.. = (\mu_n^T R_{n-1} R_{n-2} \ldots..R_0)/(\alpha_n\ldots..\alpha_1)$

or, $$\mu_0 = (R_0^T R_1^T\ldots\ldots..R_{n-1}^T\mu_n)/(\alpha_1\alpha_2\ldots\ldots..\alpha_n) \tag{106}$$

Cor 3: $\mu_{n+1}^T = \alpha_{n+1}\mu_n^T R_n^{-1} = \alpha_{n+1}(\alpha_n \quad \beta_n \quad 1)\begin{pmatrix} 0 & 1 & 0 \\ 0 & 0 & 1 \\ 1 & -a_n & -b_n \end{pmatrix}$

**Theorem 4.** If the BCF expansion of $(\alpha, \beta)$ is eventually periodic then $\alpha$ satisfies a polynomial of degree at most three[5].

Proof: Let $(\alpha, \beta) = [\{a_0, a_1, \ldots\ldots, a_n, a_{n+1}, a_{n+2}, \ldots\ldots, a_{n+m}, a_{n+1}, a_{n+2}, \ldots\ldots, a_{n+m}, \ldots\ldots\}, \{b_0, b_1, \ldots\ldots b_n, b_{n+1}, b_{n+2}, \ldots\ldots b_{n+m}, b_{n+1}, b_{n+2}, \ldots\ldots b_{n+m}, \ldots\ldots\}]$ (107)

Then, $(\alpha_{n+1}, \beta_{n+1}) = [\{a_{n+1}, a_{n+2}, \ldots\ldots, a_{n+m}, a_{n+1}, a_{n+2}, \ldots\ldots, a_{n+m}, \ldots\ldots\}, \{b_{n+1}, b_{n+2}, \ldots\ldots b_{n+m}, b_{n+1}, b_{n+2}, \ldots\ldots b_{n+m}, \ldots\ldots\}]$ (108)

By Cor 2 of Lemma 12, we get,

$$\mu_0^T = (\mu_{n+1}^T R_n R_{n-1} \ldots\ldots R_0)/(\alpha_{n+1} \alpha_n \ldots\alpha_1)$$

$$\mu_{n+1}^T = (\mu_{n+1}^T R_{n+m} R_{n+m-1} \ldots\ldots R_{n+2} R_{n+1})/(\alpha_{n+1} \alpha_{n+m}\alpha_{n+m-1} \ldots\alpha_{n+2}) \quad (109)$$

Hence,

$$\mu_0^T = (\mu_0^T R_0^{-1} R_1^{-1} \ldots R_n^{-1} R_{n+m} R_{n+m-1} \ldots R_{n+2} R_{n+1} R_n R_{n-1} \ldots R_1 R_0)/(\alpha_{n+1}\ldots\alpha_{n+m}) \quad (110)$$

This equation may be written as

$$\lambda (\alpha \quad \beta \quad 1) = (\alpha \quad \beta \quad 1) \begin{pmatrix} M_{11} & M_{12} & M_{13} \\ M_{21} & M_{22} & M_{23} \\ M_{31} & M_{32} & M_{33} \end{pmatrix} \quad (111)$$

Hence,

$$\begin{aligned} M_{11} \alpha + M_{21} \beta + M_{31} &= \lambda \alpha \\ M_{12} \alpha + M_{22} \beta + M_{32} &= \lambda \beta \\ M_{13} \alpha + M_{23} \beta + M_{33} &= \lambda \end{aligned} \quad (112)$$

Thus,

$$M_{11} \alpha + M_{21} \beta + M_{31} = (M_{13} \alpha + M_{23} \beta + M_{33}) \alpha \quad (113)$$

which implies

$$\beta = \{M_{31} + (M_{11} - M_{33})\alpha - M_{13} \alpha^2\}/(M_{23} \alpha - M_{21}) \quad (114)$$

Also,

$$M_{12} \alpha + M_{22} \beta + M_{32} = (M_{13} \alpha + M_{23} \beta + M_{33}) \beta \quad (115)$$

so that

$$M_{32} + M_{12}\,\alpha + (M_{22} - M_{33} - M_{13}\,\alpha)\{M_{31} + (M_{11} - M_{33})\alpha - M_{13}\,\alpha^2\}/(M_{23}\,\alpha - M_{21})$$
$$- M_{23}\,\{M_{31} + (M_{11} - M_{33})\,\alpha - M_{13}\,\alpha^2\}^2\,/(M_{23}\,\alpha - M_{21})^2 = 0 \tag{116}$$

Therefore,

$$(M_{32} + M_{12}\,\alpha)(M_{23}\,\alpha - M_{21})^2$$
$$+ (M_{22} - M_{33} - M_{13}\,\alpha)\{M_{31} + (M_{11} - M_{33})\,\alpha - M_{13}\,\alpha\}(M_{23}\,\alpha - M_{21})$$
$$- M_{23}\{M_{31} + (M_{11} - M_{33})\,\alpha - M_{13}\,\alpha^2\}^2 = 0 \tag{117}$$

The coefficient of $\alpha^4$ vanishes. Therefore $\alpha$ satisfies a polynomial of degree at most three.

**Conjecture**: If $\alpha$ is a cubic irrational, there exists a cubic irrational $\beta$ such that the BCF of $(\alpha, \beta)$ is eventually periodic.